\def\A{\leavevmode\setbox0\hbox{A}\lower1.4ex\hbox to\wd0{\hss`}\kern-.9\wd0A}
\def\E{\leavevmode\setbox0\hbox{E}\lower1.4ex\hbox to\wd0{\hss`\/}\kern-.9\wd0E}
\def\a{\leavevmode\setbox0\hbox{a}\lower1.4ex\hbox to\wd0{\hss`\/}\kern-\wd0a}
\def\e{\leavevmode\setbox0\hbox{e}\lower1.4ex\hbox to\wd0{\hss`\/}\kern-\wd0e}
\newcommand{\be}{\begin{equation}}
\newcommand{\ee}{\end{equation}}
\newcommand{\ba}{\begin{array}}
\newcommand{\ea}{\end{array}}
\newcommand{\beqn}{\begin{eqnarray}}
\newcommand{\eeqn}{\end{eqnarray}}

\newcommand{\zero}{\setcounter{equation}{0} \par}

\def\znakr{\raise1.5pt\hbox{\symb\char66\kern-2pt\char74}}
\def\znakl{\raise1.5pt\hbox{\symb\char73\kern-2pt\char67}}

\def\normalsize{
\setlength{\textheight}{23cm}
\setlength{\textwidth}{15cm}
\setlength{\topmargin}{-2.0cm}
\setlength{\hoffset}{-0.5cm}
\setlength{\leftmargin}{-1cm}
\setlength{\rightmargin}{2.0cm}}

\documentstyle[12pt]{article}
\normalsize

\begin{document}
\title{The Hopf algebra isomorphism between $\kappa -$Poincar\'e algebra
in the case $g^{00}=0$
and "null plane" quantum Poincar\'e algebra}
\author{ Karol Przanowski\thanks{Supported by
\L{}\'od\'z University grant no 487} \\
Department of Field Theory \\
University of \L \'od\'z \\
ul.Pomorska 149/153, 90-236 \L \'od\'z , Poland}
\date{}
\maketitle
\setcounter{section}{0}
\setcounter{page}{1}
\begin{abstract}
The Hopf algebra isomorphism between $\kappa -$Poincar\'e algebra
defined by P.Kosi\'nski and P.Ma\'slanka in the case
$g^{00}=0$~\cite{defwayla}
and "null plane" quantum Poincar\'e algebra by
A.Ballesteros, F.J.Herranz and M.A.del Olmo~\cite{null} is defined.
\end{abstract}

\section{Introduction}
\zero

Recently, considerable interest has been paid to the deformations of group
and algebras of space-time symmetries~\cite{w1}. An interesting deformation of
the Poincar\`e algebra~\cite{w2},~\cite{defwayla} as well
as group~\cite{zak} has been introduced which
depend on the dimensional deformation parameter $\kappa $; the relevant
objects are called  $\kappa -$Poincar\`e algebra and $\kappa -$Poincar\`e
group, respectively. Their structure was studied in some detail and many
of their properties are now well understood.

In the paper~\cite{null} using the so called deformation embeding method
A.Ballesteros, F.J.Herranz and M.A.del Olmo obtained the 
"null plane" quantum Poincar\'e algebra. On the other hand in~\cite{defwayla}
P.Kosi\'nski and P.Ma\'slanka presented the method for obtaining the
$\kappa -$Poincar\`e algebra for arbitrary metric tensor. It is interesting
whether the "null plane" quantum Poincar\'e algebra is the particual case
(for the special choice of the matric tensor $g^{\mu \nu}$)
of the $\kappa -$deformation presented in~\cite{defwayla}. In our paper
we solve this problem.

In the section \ref{r1} we rewrite $\kappa -$Poincar\'e algebra in the new basis.
In the section \ref{r2} we describe
the "null plane" quantum Poincar\'e algebra~\cite{null}
and define isomorpfism between these two algebras. At least
in the section \ref{r3} we define isomorphism between $\kappa -$Poincar\'e
algebra defined below and the "null plane" quantum Poincar\'e algebra.

Let us remind the definition of $\kappa -$Poincar\`e algebra .
The $\kappa -$Poincar\`e algebra $\tilde {\cal P}_\kappa $~\cite{defwayla}
(in the Majid and Ruegg basis~\cite{w9}) is a quantized universal envoloping
algebra in the sense of Drinfeld~\cite{drin} described by the
following relations:

\noindent  The commutation rules:
\beqn
\lbrack  M^{ij},P_0 \rbrack  &=& 0 ,\nonumber \\
\lbrack M^{ij},P_k  \rbrack  &=& i\kappa (\delta ^j_{\ k}
g^{0i}-\delta ^i_{\ k}g^{0j})(1-e^{-{P_0\over \kappa }})
+ i(\delta ^j_{\ k}g^{is}-\delta ^i_{\ k}g^{js})P_s ,\nonumber \\
\lbrack M^{i0},P_0  \rbrack  &=& i\kappa g^{i0}
(1-e^{-{p_0\over \kappa }})+ig^{ik}P_k ,\nonumber \\
\lbrack M^{i0},P_k  \rbrack  &=& -i{\kappa \over 2}g^{00}
\delta ^i_{\ k}(1-e^{-2{P_0\over \kappa }})-
i\delta ^i_{\ k}g^{0s}P_s e^{-{P_0\over \kappa }}+ ,\nonumber \\
 &\ & +ig^{0i}P_k (e^{-{P_0\over \kappa }}-1)+
{i\over 2\kappa }\delta ^i_{\ k}g^{rs}P_r P_s -
{i\over \kappa }g^{is} P_s P_k ,\nonumber \\
\lbrack P_\mu ,P_\nu   \rbrack  &=& 0 ,\nonumber \\
\lbrack M^{\mu \nu },M^{\lambda \sigma }  \rbrack
&=& i(g^{\mu \sigma }M^{\nu \lambda }-g^{\nu \sigma }M^{\mu \lambda }
+g^{\nu \lambda }M^{\mu \sigma }-g^{\mu \lambda }M^{\nu \sigma }) .
\label{stk}
\eeqn
The coproducts, counit and antipode:
\beqn
\Delta P_0 &=& I\otimes P_0 + P_0 \otimes I ,\nonumber \\
\Delta P_k &=& P_k \otimes e^{-{P_0\over \kappa }}+
I\otimes P_k ,\nonumber \\
\Delta M^{ij} &=& M^{ij}\otimes I + I\otimes M^{ij} ,\nonumber \\
\Delta M^{i0} &=& I\otimes M^{i0} + M^{i0}\otimes
e^{-{P_0\over \kappa }}- {1\over \kappa }M^{ij}\otimes P_j ,\nonumber \\
\varepsilon (M^{\mu \nu }) &=& 0; \ \ \
\varepsilon (P_\nu ) = 0, \nonumber \\
S(P_0) &=& -P_0 ,\nonumber \\
S(P_i) &=& -e^{P_0\over \kappa }P_i ,\nonumber \\
S(M^{ij}) &=& -M^{ij} ,\nonumber \\
S(M^{i0}) &=& -(M^{i0}+{1\over \kappa }M^{ij}P_j)e^{P_0\over \kappa },
\label{stp}
\eeqn
where $i,j,k = 1,2,3$ and
the metric tensor $g_{\mu \nu }, (\mu \nu =0,1,...,3)$ is
represented by an arbitrary nondegenerate symetric $4\times 4$ matrix
(not necessery diagonal).

\section{The $\kappa -$Poincar\'e algebra in the new basis \label{r1}}
\zero

Note that in our paper we mark
$$
A^i B^i = \sum_{n=1}^3 A^n B^n,\ A^i B_i = \sum_{n=1}^3 A^n B_n
$$
for any tensors $A^\mu ,\ B^\nu $ and $\varepsilon ^{123} = -1$.

We put:
\beqn
M^i &=& {1\over 2}\varepsilon ^{ijk}M^{jk},\
({\rm we\ have\ also\ } M^{ij}=\varepsilon ^{ijk}M^k ),  \nonumber  \\
N^i &=& M^{i0}. \nonumber
\eeqn
We define the isomorphism on generators of
$\kappa -$Poincar\'e algebra~\cite{w9}:
\beqn
{\cal P}_0 &=&  -P_0 ,\nonumber  \\
{\cal P}_i &=&  -P_i e^{P_0\over 2\kappa } ,\nonumber  \\
{\cal M}^i &=&  M^i ,\nonumber  \\
{\cal N}^i &=&  (N^i - {1\over 2\kappa }\varepsilon ^{ijk}M^j P_k)
e^{P_0\over 2\kappa }\nonumber  \\
&=& N^i e^{P_0\over 2\kappa } + {1\over 2\kappa }
\varepsilon ^{ijk}{\cal M}^j {\cal P}_k .\nonumber
\eeqn
After some calculation we get the following relations defining the 
$\kappa -$Poincar\'e algebra in the new basis:

\noindent  The commutation rules:
\newcommand{\sh}{\sinh ({ {\cal P}_0 \over 2\kappa }) }
\newcommand{\ch}{\cosh ({ {\cal P}_0 \over 2\kappa }) }
\beqn
\lbrack {\cal P}_\mu , {\cal P}_\nu \rbrack &=& 0 ,\nonumber  \\
\lbrack {\cal M}^i , {\cal P}_0 \rbrack &=& 0 ,\nonumber  \\
\lbrack {\cal M}^i , {\cal P}_k \rbrack &=& i\varepsilon ^{ijl}\delta ^l_k
(2\kappa g^{oj}\sh +g^{js}{\cal P}_s) ,\nonumber  \\
\lbrack {\cal N}^i , {\cal P}_0 \rbrack &=&
2i\kappa g^{i0}\sh + ig^{ik}{\cal P}_k ,\nonumber  \\
\lbrack {\cal N}^i , {\cal P}_k \rbrack &=& -i\kappa g^{00}\delta ^i_k
\sinh({{\cal P}_0\over \kappa }) -
i\delta ^i_k g^{0s}{\cal P}_s \ch ,\nonumber  \\
\lbrack {\cal M}^i , {\cal M}^j \rbrack &=& -
i\varepsilon ^{ijk}g^{ks}{\cal M}^s ,\nonumber  \\
\lbrack {\cal N}^i , {\cal M}^j \rbrack &=&
i\varepsilon ^{jrs}g^{ir}{\cal N}^s
+ig^{i0}{\cal M}^j \ch -i\delta ^{ij}g^{k0}{\cal M}^k \ch ,\nonumber  \\
\lbrack {\cal N}^i , {\cal N}^j \rbrack &=&
ig^{j0}{\cal N}^i \ch -ig^{i0}{\cal N}^j \ch
-{i\over 4\kappa ^2}\varepsilon ^{ijs}
g^{kr}{\cal M}^r {\cal P}_s {\cal P}_k \nonumber \\
&\ &+{i\over 2\kappa }\varepsilon ^{jrs}g^{i0}{\cal M}^r {\cal P}_s \sh
-{i\over 2\kappa }\varepsilon ^{irs}
g^{j0}{\cal M}^r {\cal P}_s \sh   \nonumber  \\
&\ &-i\varepsilon ^{ijk}g^{00}{\cal M}^k
\cosh ({{\cal {\cal P}}_0\over \kappa })
-{i\over \kappa }\varepsilon ^{ijk}
{\cal M}^k g^{s0}{\cal P}_s \sh  .
\label{newk}
\eeqn
The coproducts, counit and antipode:
\beqn
\Delta {\cal P}_0 &=& {\cal P}_0 \otimes I + I\otimes {\cal P}_0 ,
\nonumber  \\
\Delta {\cal P}_i &=& {\cal P}_i \otimes e^{{\cal P}_0\over 2\kappa } +
e^{-{{\cal P}_0\over 2\kappa }}\otimes {\cal P}_i ,\nonumber  \\
\Delta {\cal M}^i &=& {\cal M}^i \otimes I + I\otimes {\cal M}^i ,
\nonumber  \\
\Delta {\cal N}^i &=& {\cal N}^i \otimes e^{{{\cal P}_0\over 2\kappa }} +
e^{-{{\cal P}_0\over 2\kappa }}\otimes {\cal N}^i -
{1\over 2\kappa }\varepsilon ^{ijk}(e^{-{{\cal P}_0\over 2\kappa }}
{\cal M}^j \otimes {\cal P}_k - {\cal P}_k \otimes {\cal M}^j
e^{{{\cal P}_0\over 2\kappa }}),\nonumber  \\
\varepsilon ({\cal X}) &=& 0,\ {\rm for}\
{\cal X}={\cal N}^i,\ {\cal M}^i,\ {\cal P}_\mu , \nonumber  \\
S({\cal X}) &=& -{\cal X},\ {\rm for}\
{\cal X}={\cal M}^i,\ {\cal P}_\mu , \nonumber  \\
S({\cal N}^i) &=& -{\cal N}^i +3i(g^{i0}\sh +{1\over 2\kappa }
g^{ik}{\cal P}_k ), \nonumber \\
({\rm or}\ S({\cal X}) &=& - e^{{3{\cal P}_0\over 2\kappa }}
{\cal X} e^{-{3{\cal P}_0\over 2\kappa }},\ {\rm for}\
{\cal X}= {\cal M}^i,\ {\cal N}^i,\ {\cal P}_\mu ).
\label{newp}
\eeqn

Note, if we take diagonal metric tensor $g_{\mu \nu }={\rm Diag}(1,-1,-1,-1)$,
we obtain $\kappa -$Poincar\'e algebra considered
in~\cite{w9},~\cite{w2},~\cite{mas}.

\section{The "null plane" quantum Poincar\'e algebra and isomorphism
\label{r2}}
\zero

The "null plane" quantum Poincar\'e algebra is a Hopf *-algebra generated
by ten elements:
$P_+ ,\ P_- ,\ P_1 ,\ P_2 ,\ E_1 ,\ E_2 ,\ F_1 ,\ F_2 ,\ J_3 , K_3$
and $z$-deformation parameter
with the following relations~\cite{null}:

\noindent  The commutation rules:
\newcommand{\shn}{\sinh (zP_+) }
\newcommand{\chn}{\cosh (zP_+) }
\beqn
\lbrack K_3 , P_+ \rbrack &=& {\shn \over z} ,\nonumber  \\
\lbrack K_3 , P_- \rbrack &=& -P_- \chn ,\nonumber  \\
\lbrack K_3 , E_i \rbrack &=& E_i \chn ,\nonumber  \\
\lbrack K_3 , F_1 \rbrack &=& -F_1 \chn + zE_1 P_- \shn
- z^2 P_2 W^q_+ ,\nonumber  \\
\lbrack K_3 , F_2 \rbrack &=& -F_2 \chn + zE_2 P_- \shn
- z^2 P_1 W^q_+ ,\nonumber  \\
\lbrack J_3 , P_i \rbrack &=& -\varepsilon _{ij3}P_j ,\nonumber  \\
\lbrack J_3 , E_i \rbrack &=& -\varepsilon _{ij3}E_j ,\nonumber  \\
\lbrack J_3 , F_i \rbrack &=& -\varepsilon _{ij3}F_j ,\nonumber  \\
\lbrack E_i , P_j \rbrack &=& \delta _{ij}{\shn \over z} ,\nonumber  \\
\lbrack F_i , P_j \rbrack &=& \delta _{ij}P_- \chn ,\nonumber  \\
\lbrack E_i , F_j \rbrack &=& \delta _{ij}K_3
+ \varepsilon _{ij3}\chn ,\nonumber  \\
\lbrack P_+ , F_i \rbrack &=& -P_i ,\nonumber  \\
\lbrack F_1 , F_2 \rbrack &=&  z^2 P_- W^q_+ + zP_- J_3 \shn ,\nonumber  \\
\lbrack P_- , E_i \rbrack &=& -P_i .
\label{nullk}
\eeqn
The coproducts, counit and antipode:
\beqn
\Delta X &=& I\otimes X + X\otimes I,\ {\rm for}
\ X=P_+,\ E_i,\ J_3 ,\nonumber  \\
\Delta Y &=& e^{-zP_+}\otimes Y + Y\otimes e^{zP_+},\ {\rm for}
\ Y=P_-,\ P_i ,\nonumber  \\
\Delta F_1 &=& e^{-zP_+}\otimes F_1 + F_1 \otimes e^{zP_+}
+ ze^{-zP_+}E_1 \otimes P_-  \nonumber  \\
&\ & - zP_- \otimes E_1 e^{zP_+}
+ ze^{-zP_+}J_3 \otimes P_2 - zP_2 \otimes J_3 e^{zP_+} ,\nonumber  \\
\Delta F_2 &=& e^{-zP_+}\otimes F_2 + F_2 \otimes e^{zP_+}
+ ze^{-zP_+}E_2 \otimes P_-  \nonumber  \\
&\ & - zP_- \otimes E_2 e^{zP_+}
- ze^{-zP_+}J_3 \otimes P_1 + zP_1 \otimes J_3 e^{zP_+} ,\nonumber  \\
\Delta K_3 &=& e^{-zP_+}\otimes K_3 + K_3 \otimes e^{zP_+}
+ ze^{-zP_+}E_1 \otimes P_1  \nonumber  \\
&\ & - zP_1 \otimes E_1 e^{zP_+}
+ ze^{-zP_+}E_2 \otimes P_2 - zP_2 \otimes E_2 e^{zP_+} ,\nonumber  \\
\varepsilon (X) &=& 0,\ S(X)= -e^{3zP_+}Xe^{-3zP_+},\ {\rm for}
\ X=P_\pm ,\ P_i ,\ F_i ,\ E_i , \ J_3 , \ K_3 ,
\label{nullp}
\eeqn
where $W^q_+ = E_1 P_2 - E_2 P_1 + J_3 {\shn \over z}$ and $i,j = 1,2$.

If we subtitude the followig expresions into
the relations (\ref{newk}),(\ref{newp}) describing
the $\kappa -$Poincar\'e algebra:
\beqn
g^{\mu \nu } &=& \left\lbrack
\ba{rrrr}
0 & 0 & 0 & 1 \\
0 & -1 & 0 & 0 \\
0 & 0 & -1 & 0 \\
1 & 0 & 0 & 0
\ea
\right\rbrack
\label{gmunu}
\eeqn
and:
\beqn
P_+ &=& {\cal P}_0 ,\nonumber \\
P_- &=& {\cal P}_3 ,\nonumber \\
P_i &=& {\cal P}_i ,\nonumber \\
K_3 &=& -i{\cal N}^3 ,\nonumber \\
F_i &=& i{\cal N}^i ,\nonumber \\
E_1 &=& i{\cal M}^2 ,\nonumber \\
E_2 &=& -i{\cal M}^1 ,\nonumber \\
J_3 &=& -i{\cal M}^3 ,\nonumber \\
z &=& {1\over 2\kappa } ,\nonumber
\eeqn
we get the relations
(\ref{nullk}),(\ref{nullp}) describing the
"null plane" quantum Poincar\'e algebra.

\section{Summary \label{r3}}
\zero

We define Hopf algebra isomorpfism from $\kappa -$Poincar\'e algebra
(\ref{stk}),(\ref{stp}) for metric tensor (\ref{gmunu}) to
the "null plane" quantum Poincar\'e algebra by putting:
\beqn
P_+ &=& -P_0 ,\nonumber \\
P_- &=& -P_3 e^{P_0\over 2\kappa } ,\nonumber \\
P_i &=& -P_i e^{P_0\over 2\kappa } ,\nonumber \\
K_3 &=& -i(M^{30} + {1\over 2\kappa }M^{3k}P_k)
e^{P_0\over 2\kappa } ,\nonumber \\
F_i &=& i(M^{i0} + {1\over 2\kappa }M^{ik}P_k)
e^{P_0\over 2\kappa }  ,\nonumber \\
E_1 &=& {i\over 2}\varepsilon ^{2jk}M^{jk} ,\nonumber \\
E_2 &=& -{i\over 2}\varepsilon ^{1jk}M^{jk} ,\nonumber \\
J_3 &=& -{i\over 2}\varepsilon ^{3jk}M^{jk} ,\nonumber \\
z &=& {1\over 2\kappa } ,\nonumber
\eeqn
for $i=1,2;\ j,k=1,2,3$.

We conclude that the "null plane" quantum Poincar\'e algebra is
a case of the $\kappa -$Poincar\'e algebra for $g^{00}=0$.

\section{Acknowledgment}
\zero
I would like to thank prof.Jerzy Lukierski for suggesting that I
undertake this problem.
I would also like to ackoweledge that prof.Jerzy Lukierski and
dr Pawe\l{}  Ma\`slanka kindly discussed this problem with me.


\end{document}